\documentclass[letterpaper, 10 pt, conference]{ieeeconf}  
\IEEEoverridecommandlockouts                              

\usepackage{graphics} 
\usepackage{epsfig} 
\usepackage{times} 
\usepackage{mathtools}
\usepackage{amsmath} 
\usepackage{amssymb}  
\usepackage{amsthm}
\usepackage{color}

\usepackage{tikz}
\usetikzlibrary{shapes,arrows}
\usepackage{pgfplots} 
\usepackage{pgfgantt}
\usepackage{pdflscape}
\usepackage{nicefrac}
\usepackage{prettyref}
\pgfplotsset{compat=1.5}
\pgfplotsset{plot coordinates/math parser=false}
\newlength\fwidth
\usepackage{grffile}
\usetikzlibrary{plotmarks}
\usetikzlibrary{arrows.meta}
\usepgfplotslibrary{patchplots}
\usepackage{units}
\let\labelindent\relax
\usepackage{enumitem}
\usepackage{hyperref}
\usepackage{float}
\usepackage{subcaption}
\usepackage{algorithm}
\usepackage{algpseudocode}
\usepackage{cite}
\newtheorem{definition}{Definition}
\newtheorem{remark}{Remark}

\newtheorem{problem}{Problem}
\newtheorem{proposition}{Proposition}
\newtheorem{assumption}{Assumption}

\newcommand{\defeq}{\vcentcolon=}
\newcommand{\T}{\scriptscriptstyle\top}       

\newcommand{\mathst}{\text{s.t.}}

\DeclareMathOperator*{\argmax}{arg\,\operatorname*{max}}

\DeclareUnicodeCharacter{2212}{-}

\algnewcommand{\algorithmicgoto}{\textbf{go to}}%
\algnewcommand{\Goto}[1]{\algorithmicgoto~\ref{#1}}

\makeatletter
\def\widebreve{\mathpalette\wide@breve}
\def\wide@breve#1#2{\sbox\z@{$#1#2$}%
     \mathop{\vbox{\m@th\ialign{##\crcr
\kern0.18em\vspace{-0.016cm}\brevefill#1{0.55\wd\z@}\crcr\noalign{\nointerlineskip}%
                    $\hss#1#2\hss$\crcr}}}\limits}
\def\brevefill#1#2{$\m@th\sbox\tw@{$#1($}%
  \hss\resizebox{#2}{\wd\tw@}{\rotatebox[origin=c]{90}{\upshape(}}\hss$}
\makeatletter                                                          

\IEEEoverridecommandlockouts                              
\title{\LARGE \bf
Optimization-based Verification of Discrete-time Control Barrier Functions: A Branch-and-Bound Approach
}

\author{Erfan Shakhesi, W.P.M.H. (Maurice) Heemels, and Alexander Katriniok
\thanks{The research has received funding from the European Research Council under the Advanced ERC Grant Agreement PROACTHIS, no. 101055384.}
\thanks{The authors are with the Control Systems Technology section, Mechanical Engineering, Eindhoven
University of Technology, The Netherlands. E-mail: {\{\tt\small e.shakhesi, m.heemels, a.katriniok\}@tue.nl}}
}

\usepackage{eso-pic}
\AddToShipoutPictureBG*{%
\AtPageUpperLeft{%
\setlength\unitlength{1in}%
\hspace*{\dimexpr0.5\paperwidth\relax}
\makebox(0,-1.05)[c]{
\begin{tabular}{c c}

To appear in the Proceedings of the 63rd IEEE Conference on Decision and Control, Italy, 2024. \\
Uploaded to arXiv on September 23, 2024.

\end{tabular}}}}

\AddToShipoutPictureBG*{%
\AtPageUpperLeft{%
\setlength\unitlength{1in}%
\hspace*{\dimexpr0.5\paperwidth\relax}
\makebox(0,-21.25)[c]{
\footnotesize
\begin{tabular}{c c}
© 2024 IEEE. Personal use of this material is permitted. Permission from
IEEE must be obtained for all other uses, in any \\
current or future media, including reprinting/republishing
this material for advertising or promotional purposes, creating new \\
collective works, for resale or redistribution to servers or lists,
or reuse of any copyrighted component of this work in other works.
\end{tabular}}}}

\begin{document}

\maketitle
\thispagestyle{empty}
\pagestyle{empty}

\begin{abstract}
Discrete-time Control Barrier Functions (DTCBFs) form a powerful control theoretic tool to guarantee safety and synthesize safe controllers for discrete-time dynamical systems. In this paper, we provide an optimization-based algorithm, inspired by the $\alpha$BB algorithm, for the verification of a candidate DTCBF, i.e., either verifying a given candidate function as a valid DTCBF or falsifying it by providing a counterexample for a general nonlinear discrete-time system with input constraints. This method is applicable whether a corresponding control policy is known or unknown. We apply our method to a numerical case study to illustrate its efficacy.
\end{abstract}

\section{Introduction}
\label{sec:sec1}
Safety-critical systems are characterized as those in which a failure could result in substantial harm or damage. Such systems can be identified across diverse domains, including aerospace, healthcare, and automotive. Due to the severe consequences of safety violations, the availability of design methods for controllers with formal safety guarantees for safety-critical systems is crucial. 

Typically, dynamical systems are considered safe, if their system trajectories remain in predefined safe sets. Thus, guaranteeing safety can be achieved by constructing controlled invariant sets within these safe sets. To achieve this and to synthesize safe controllers, Control Barrier Functions (CBFs) are widely recognized in the literature as a powerful tool \cite{Ames2014a}. In particular, CBF-based controllers often act as safety filters, adjusting nominal control inputs to prevent the system state from leaving a safe set. This control technique was initially developed for continuous-time systems \cite{Ames2014a}, but an extension of CBFs to discrete-time systems has been introduced in \cite{Agrawal2017a}, which is referred to as discrete-time CBFs (DTCBFs). DTCBFs are also used within Model Predictive Control (MPC) to formulate stage constraints \cite{Zeng2021a}, \cite{Zeng2021b}, but also to guarantee recursive feasibility \cite{katriniok2023}.

Synthesizing DTCBFs for arbitrary discrete-time dynamical systems with input constraints presents inherent challenges and remains unsolved in the literature. As a result, almost all the previous works that utilized DTCBFs, including \cite{Agrawal2017a}, \cite{Zeng2021a}, and \cite{Zeng2021b}, either operated under the assumption that DTCBFs are given and used them to synthesize safe controllers or employed \textit{candidate} DTCBFs that do not necessarily adhere to the formal DTCBF definition rigorously, thereby lacking guarantees of safety and feasibility. Candidate DTCBFs can be constructed, for example, through handcrafting \cite{katriniok2023} or by employing learning-based methods \cite{Didier2023}. To formally guarantee safety in these cases, an algorithm is required to verify whether the candidate DTCBF is valid by adhering to the formal DTCBF definition. 

For continuous-time CBFs, \cite{Zhang2023, Clark2021, Clark2022, Isaly2022} propose methods using Sum-of-Squares (SOS) programming for the verification of polynomial candidate CBFs considering polynomial system dynamics. Additionally and closer to this paper is the optimization-based method using grid sampling proposed in \cite{Tan2022} for the verification of multiple candidate continuous-time CBFs considering control-affine systems. This method computes a control input that satisfies the CBF constraint at a fixed sample state by solving an optimization problem, which is convex under the assumption that the control admissible set is convex. Then, by leveraging the Lipschitz continuity of the system dynamics, it determines a region in the state space where this control input remains valid. This process is repeated at different sample states until all states are covered. However, for the verification of DTCBFs, the resulting optimization problem at a fixed state is generally non-convex, 
rendering the method in \cite{Tan2022} inapplicable.
To the best of our knowledge, there is currently no systematic method available in the literature for the verification of candidate DTCBFs. Nevertheless, if a corresponding control policy is provided, nonlinear optimization solvers, such as the $\alpha$BB ($\alpha$-based branch-and-bound) algorithm \cite{Floudas2010a} or SMT (Satisfiability Modulo Theories) solvers, including dReal \cite{dReal}, can be utilized, although potentially with some limitations, as discussed in this paper. Moreover, ideas from the verification of a candidate discrete-time Control Lyapunov Function (CLF) with a given control policy, as in \cite{wu2023neural}, could be used, since these are closely related to the verification of a candidate DTCBF with a known control policy. 
However, in real-world scenarios where a corresponding control policy is often unknown, these approaches cannot be directly used to verify candidate DTCBFs.

In this paper, we address the above mentioned challenge, by providing the following key contributions: 
\begin{enumerate}[leftmargin=15pt]
    \item We introduce a novel branch-and-bound (BB) method, inspired by the $\alpha$BB algorithm \cite{Floudas2010a}, to either verify a candidate DTCBF with a given control policy as a valid DTCBF for a dynamical system with input constraints, or falsify it by providing a counterexample.
    \item 
    We extend our BB method to the case where the control policy is unknown. In this case, a corresponding piecewise constant control policy that satisfies input constraints is obtained upon verifying a candidate DTCBF.
    \item We apply our method to a numerical case study to verify or falsify a candidate DTCBF, which cannot be done using existing works in the literature.
\end{enumerate}

\textit{Notation:} We use $\mathbb{R}$, $\mathbb{R}_{>0}$, and $\mathbb{R}_{\geqslant 0}$ to denote the set of real numbers, positive real numbers, and non-negative real numbers, respectively. Additionally, $\mathbb{R}^{n}$ represents the set of all $n$-dimensional vectors of real numbers. For a vector $x \in \mathbb{R}^n$, $x_i \in \mathbb{R}$, $i \in \{1, \hdots, n\}$, represents the \mbox{$i$-th} element of $x$, and $x^{\T}$ represents its transpose. Moreover, $\mathbb{N} \defeq \{1,2,3, \hdots\}$, and $\mathbb{N}_0 \defeq \mathbb{N} \cup \{0\}$. 
\begin{definition}[Zero-superlevel Set]
    The zero-superlevel set $\mathcal{C}$ of a function $h:\mathbb{R}^n \rightarrow \mathbb{R}$ is defined as 
    \begin{align}
        \mathcal{C} \defeq \{ x \in \mathbb{R}^n \mid h(x) \geqslant 0\}. \nonumber
    \end{align}
\end{definition}
\begin{definition}[Class-$\mathcal{K}_\infty$ Functions]
    A continuous function $\gamma:\mathbb{R}_{\geqslant 0} \rightarrow \mathbb{R}_{\geqslant 0}$ is said to be a \textit{class-$\mathcal{K}_\infty$} function, denoted by $\gamma \in \mathcal{K}_{\infty}$, if it is strictly increasing, $\gamma(0) = 0$, and $\lim_{r\rightarrow \infty} \gamma(r)=\infty$.
\end{definition}
For $\gamma \in \mathcal{K}_{\infty}$, we utilize the notation $\gamma \in \mathcal{K}^{\leqslant \mathrm{id}}_{\infty}$ to indicate that $\gamma(r) \leqslant r$ for all $r \in \mathbb{R}_{\geqslant 0}$.
\begin{definition}[$n$-rectangle Set]
    We define $X \subset \mathbb{R}^n$ as an $n$-rectangle set, utilizing the notation 
        $X \defeq [x^{lb}, ~ x^{ub}],$ 
    where $x^{lb},x^{ub} \in \mathbb{R}^n$ with $x_i^{lb} \leqslant x_i^{ub}$, $i \in \{1, \hdots, n\}$, given by 
    \begin{align}
        X \defeq \{ x \in \mathbb{R}^n &\mid x_i^{lb} \leqslant x_i \leqslant x_i^{ub}, ~i \in \{1, \hdots, n\}\}. \nonumber
    \end{align}
\end{definition}

\section{Background on Control Barrier Functions and Problem Statement} \label{sec:sec2}
We consider discrete-time systems of the form
\begin{align} \label{eq:sec2:discrete-time-dynamical-system}
    x^+ = f(x, u),
\end{align}
with state vector $x \in \mathbb{R}^n$, control input vector \mbox{$u \in \mathbb{U} \subseteq \mathbb{R}^m$}, both at the current time instant, state vector $x^+ \in \mathbb{R}^n$ at the next time instant, and mapping \mbox{$f:\mathbb{R}^n\times \mathbb{R}^m \rightarrow \mathbb{R}^n$}. Here, $\mathbb{U} \subseteq \mathbb{R}^m$ is the control admissible set. 
\begin{definition}[Controlled Invariance \cite{Ames2019a}]
    For the system \eqref{eq:sec2:discrete-time-dynamical-system} with the control admissible set $\mathbb{U}$, a set $\mathcal{C} \subset \mathbb{R}^n$ is controlled invariant, if, for every $x \in \mathcal{C}$, there exists a control input $u \in \mathbb{U}$ such that $f(x,u) \in \mathcal{C}$.
\end{definition}
To construct controlled invariant sets, a method that gained increasing attention in the literature is to utilize DTCBFs. 
\begin{definition}[DTCBF \cite{Agrawal2017a}, \cite{Zeng2021a}] \label{def:sec2:DTCBF}
    Consider a function \mbox{$h:\mathbb{R}^n \rightarrow \mathbb{R}$} with zero-superlevel set $\mathcal{C}$. For the system \eqref{eq:sec2:discrete-time-dynamical-system} with the control admissible set $\mathbb{U}$, $h$ is a \textit{discrete-time Control Barrier Function (DTCBF)}, if there exists a $\gamma \in \mathcal{K}^{\leqslant \mathrm{id}}_{\infty}$ such that for every $x \in \mathcal{C}$, there exists a control input $u \in \mathbb{U}$ with
    \begin{align} \label{eq:sec2:DTCBF-constraint}
        h(f(x,u)) - h(x) \geqslant - \gamma(h(x)).
    \end{align}
\end{definition}
In addition to the fact that the zero-superlevel set $\mathcal{C}$ of a DTCBF $h$ is \textit{controlled invariant} \cite{Agrawal2017a}, the DTCBF \mbox{constraint \eqref{eq:sec2:DTCBF-constraint}} incorporates an additional term, represented by $\gamma \in \mathcal{K}^{\leqslant \mathrm{id}}_{\infty}$. This term regulates the rate at which the states of the system \eqref{eq:sec2:discrete-time-dynamical-system} can approach the boundary of $\mathcal{C}$. The advantages of adjusting $\gamma$ are discussed in more detail in \cite{Zeng2021a, katriniok2023}. In this paper, we assume that $\gamma \in \mathcal{K}^{\leqslant \mathrm{id}}_{\infty}$ in \mbox{Definition \ref{def:sec2:DTCBF}} is \textit{a priori} given and fixed, similar to \cite{chen2021, tonkens2022refining, Tan2022}, and, we say that $h$ with a given $\gamma$, denoted by $(h,\gamma)$, is a DTCBF for the system \eqref{eq:sec2:discrete-time-dynamical-system} with the control admissible set $\mathbb{U}$, if for every $x \in \mathcal{C}$, there exists $u \in \mathbb{U}$ such that the DTCBF \mbox{constraint \eqref{eq:sec2:DTCBF-constraint}} is satisfied. Additionally, we say that a control policy $\pi:\mathcal{C} \rightarrow \mathbb{U}$ is a \textit{friend} of $(h, \gamma)$, assuming $\mathbb{U}$ and $f$ are clear from the context, if for all $x \in \mathcal{C}$,
    \begin{align}
        \pi(x) \in \mathbb{U} \text{~~and~~} h(f(x,\pi(x))) - h(x) \geqslant -\gamma(h(x)). \nonumber
    \end{align}

In the following, we explain the verification problem and assumptions we consider to solve this problem.
\begin{problem}[Verification] \label{problem:verification}
    Consider a candidate DTCBF \mbox{$h:\mathbb{R}^n \rightarrow \mathbb{R}$} and a function \mbox{$\gamma \in \mathcal{K}^{\leqslant \mathrm{id}}_{\infty}$}, which are assumed to be given. The objective is to either verify $(h, \gamma)$ as a valid DTCBF for the system \eqref{eq:sec2:discrete-time-dynamical-system} with the control admissible set $\mathbb{U} \subseteq \mathbb{R}^m$, or falsify it by providing a counterexample, in case a control policy $\pi:\mathcal{C} \rightarrow \mathbb{U}$ is 
    \begin{enumerate}[label=(\roman*), ref=(\roman*)]
        \item \textit{a priori} given, or \label{case:problem1:known}
        \item unknown. \label{case:problem1:unknown}
    \end{enumerate}
\end{problem}
\begin{assumption}
    The zero-superlevel set $\mathcal{C}$ of $h$, and the control admissible set $\mathbb{U}$ are compact.
\end{assumption}
\begin{assumption}
    The mapping $f$ associated with the system \eqref{eq:sec2:discrete-time-dynamical-system}, the candidate DTCBF $h$, and the control policy $\pi$, if given, are continuous.
\end{assumption}

\section{Verification of Candidate DTCBFs} \label{sec:sec4}
To solve Problem \ref{problem:verification}, we frame it as two separate optimization problems: one for \mbox{Case \ref{case:problem1:known}} and one for Case \ref{case:problem1:unknown}.
\begin{proposition}[Known Control Policy] \label{proposition:sec3:verification-known}
    Consider a candidate DTCBF $h:\mathbb{R}^n \rightarrow \mathbb{R}$, its zero-superlevel set $\mathcal{C}$, with a $\gamma \in \mathcal{K}^{\leqslant \mathrm{id}}_{\infty}$ and a control policy $\pi:\mathcal{C} \rightarrow \mathbb{U}$. Additionally, consider an $n$-rectangle set $\mathbb{X}$ such that $\mathcal{C} \subseteq \mathbb{X}$. For the system \eqref{eq:sec2:discrete-time-dynamical-system} with the control admissible set $\mathbb{U}$, $(h,\gamma)$ is a DTCBF and $\pi$ is a friend of $(h, \gamma)$, if and only if $\mathcal{F}^* \geqslant 0$, where $\mathcal{F}^* \in \mathbb{R}$ is the global minimum of
    \begin{subequations} \label{eq:sec3:prop:verification-optimization-known}
        \begin{align}
            \mathcal{F}^* \defeq \min_{x \in \mathbb{X}} ~&~ h(f(x,\pi(x))) - h(x) + \gamma(h(x)) \label{eq:sec3:prop:verification-optimization-known-objective} \\
            \mathrm{s.t.} ~&~ -h(x) \leqslant 0.  \label{eq:sec3:prop:verification-optimization-known-constraint}
        \end{align}
    \end{subequations}
    \begin{proof}
        If $\mathcal{F}^* \geqslant 0$, it implies that the DTCBF constraint \eqref{eq:sec2:DTCBF-constraint} is satisfied with the control policy $\pi$ for all $x \in \mathcal{C}$. Thus, $(h, \gamma)$ is a DTCBF according to \mbox{Definition \ref{def:sec2:DTCBF}}, and $\pi$ is a friend of $(h, \gamma)$. To prove the converse, if $(h, \gamma)$ is a DTCBF and $\pi$ is a friend of $(h, \gamma)$, then the DTCBF constraint \eqref{eq:sec2:DTCBF-constraint} is satisfied for all $x \in \mathcal{C}$. Consequently, it holds that $\mathcal{F}^* \geqslant 0$.
    \end{proof}
\end{proposition}
\begin{proposition}[Unknown Control Policy] \label{prop:sec3:verification-unknown}
    Consider a candidate DTCBF $h:\mathbb{R}^n \rightarrow \mathbb{R}$, its zero-superlevel set $\mathcal{C}$, with a \mbox{$\gamma \in \mathcal{K}^{\leqslant \mathrm{id}}_{\infty}$}. Additionally, consider an $n$-rectangle set $\mathbb{X}$ such that $\mathcal{C} \subseteq \mathbb{X}$. For the system \eqref{eq:sec2:discrete-time-dynamical-system} with the control admissible set $\mathbb{U}$, $(h, \gamma)$ is a DTCBF, if and only if $\mathcal{F}^* \geqslant 0$, where $\mathcal{F}^* \in \mathbb{R}$ is the global optimum of
    \begin{subequations} \label{eq:sec3:prop:verification-optimization-unknown}
    \begin{align}
        \mathcal{F}^* \defeq \min_{x \in \mathbb{X}} \max_{u \in \mathbb{U}} ~&~ h(f(x, u)) - h(x) + \gamma(h(x)) \\
        \mathrm{s.t.} ~&~ -h(x) \leqslant 0.
    \end{align}
    \end{subequations}
\begin{proof}
    The inner maximization problem in \eqref{eq:sec3:prop:verification-optimization-unknown} corresponds to finding a suitable control input $u \in \mathbb{U}$ that satisfies the DTCBF constraint \eqref{eq:sec2:DTCBF-constraint} for a particular $x \in \mathcal{C}$. Then, $\mathcal{F}^* \geqslant 0$ implies that for all $x \in \mathcal{C}$, there exists $u \in \mathbb{U}$ that satisfies the DTCBF constraint \eqref{eq:sec2:DTCBF-constraint}, which implies that $(h,\gamma)$ is a DTCBF. Conversely, if $(h,\gamma)$ is a DTCBF according to Definition \ref{def:sec2:DTCBF}, then for every $x \in \mathcal{C}$, there exists a control input $u \in \mathbb{U}$ that satisfies the DTCBF constraint \eqref{eq:sec2:DTCBF-constraint}. Hence, $\mathcal{F}^* \geqslant 0$. 
\end{proof}
\end{proposition}

\section{Numerical solution of the verification problem}
\label{sec:sec4b}
Since the optimization problems \eqref{eq:sec3:prop:verification-optimization-known} and \eqref{eq:sec3:prop:verification-optimization-unknown} are generally non-convex, we aim to propose a BB method to determine whether $\mathcal{F}^*$ is non-negative.
\subsection{Preliminaries}\label{sec3:preliminaries}
First, we introduce some necessary fundamentals.
\begin{definition}[Convex Underestimators \cite{ADJIMAN1998a}]
    A function \mbox{$\widebreve{\mathcal{F}}:\mathbb{R}^n \rightarrow \mathbb{R}$} is a convex underestimator of a given function \mbox{$\mathcal{F}:\mathbb{R}^n \rightarrow \mathbb{R}$} on a set $X$, if $\widebreve{\mathcal{F}}$ is convex on $X$, and \mbox{$\widebreve{\mathcal{F}}(x) \leqslant \mathcal{F}(x)$} for all $x \in X$. 
\end{definition}
In the literature, a common approach to constructing a convex underestimator of a function $\mathcal{F}:\mathbb{R}^n \rightarrow \mathbb{R}$ on an $n$-rectangle set $X \defeq [x^{lb}, ~ x^{ub}]$ involves adding a quadratic polynomial to $\mathcal{F}$ \cite{ADJIMAN1998a}. This polynomial is chosen such that it is non-positive on $X$ and has sufficiently large coefficients for the quadratic terms to overcome the non-convexity of $\mathcal{F}$. Specifically, the function $\widebreve{\mathcal{F}}:\mathbb{R}^n \rightarrow \mathbb{R}$, defined as 
\begin{align} \label{eq:sec4:preliminaries:convex-underestimator}
    \widebreve{\mathcal{F}}(x) \defeq \mathcal{F}(x) + \sum_{i=1}^{n} \alpha_i\left(x_i^{lb} - x_i\right)\left(x_i^{ub} - x_i\right),
\end{align}
is a convex underestimator of $\mathcal{F}$ on $X \defeq [x^{lb}, ~ x^{ub}]$, if $\alpha_i \in \mathbb{R}_{\geqslant 0}$, $i \in \{1, \hdots, n\}$, are sufficiently large. Numerous ways to compute suitable values of $\alpha_i$, $i \in \{1,\hdots, n\}$, can be found in \cite{ADJIMAN1998a}, such as the Scaled Gerschgorin proposition. \\
Two properties of convex underestimators, constructed according to \eqref{eq:sec4:preliminaries:convex-underestimator}, are now recalled from \cite{Maranas1994b}.
\begin{proposition}[Maximum Separation \cite{Maranas1994b}] \label{prop:sec4:preliminaries:max-separation}
    The maximum separation between a function $\mathcal{F}:\mathbb{R}^n \rightarrow \mathbb{R}$ and a convex underestimator $\widebreve{\mathcal{F}}:\mathbb{R}^n \rightarrow \mathbb{R}$, as in \eqref{eq:sec4:preliminaries:convex-underestimator} on an $n$-rectangle set $X \defeq [x^{lb}, \; x^{ub}]$, occurs at the middle of $X$, that is, it holds
    \begin{align}
        \max_{x\in X} ~ (\mathcal{F}(x) - \widebreve{\mathcal{F}}(x)) = \frac{1}{4}\sum_{i=1}^{n} \alpha_i\left(x_i^{ub} - x_i^{lb}\right)^2. \nonumber
    \end{align}
\end{proposition}
\begin{proposition}[Tighter Convex Underestimators \cite{Maranas1994b}] \label{prop:sec4:preliminaries:tigher-convex-understimator}
    Consider two $n$-rectangle sets $X^{(1)} \defeq [x^{lb, (1)}, ~ x^{ub, (1)}]$ and $X^{(2)} \defeq [x^{lb, (2)}, ~ x^{ub, (2)}]$ with \mbox{$X^{(2)} \subseteq X^{(1)}$}. Let $\widebreve{\mathcal{F}}^{(1)}$ and $\widebreve{\mathcal{F}}^{(2)}$ be convex underestimators of a function $\mathcal{F}: \mathbb{R}^n \rightarrow \mathbb{R}$ on $X^{(1)}$ and $X^{(2)}$, respectively, constructed as
    \begin{align}
        \widebreve{\mathcal{F}}^{(l)}(x) &\defeq \mathcal{F}(x) + \sum_{i=1}^{n} \alpha_i^{(l)}\bigl(x_i^{lb, (l)} - x_i\bigr)\bigl(x_i^{ub,(l)} - x_i\bigr), \nonumber 
    \end{align}
    $l \in \{1, 2\}$, where $\alpha_i^{(1)}, \alpha_i^{(2)} \in \mathbb{R}_{\geqslant 0}$ are computed based on one of the methods discussed in \cite{ADJIMAN1998a}.
    Then, 
        \mbox{$\alpha_i^{(2)} \leqslant \alpha_i^{(1)}$} for all $i \in \{1,\hdots, n\}$, 
    and $\widebreve{\mathcal{F}}^{(2)}$ is a tighter convex underestimator of $\mathcal{F}$ than $\widebreve{\mathcal{F}}^{(1)}$ on $X^{(2)}$, in the sense that, for all $x \in X^{(2)}$,
        $\widebreve{\mathcal{F}}^{(1)}(x) \leqslant \widebreve{\mathcal{F}}^{(2)}(x) \leqslant \mathcal{F}(x)$. 
\end{proposition}
\subsection{Proposed Verification Algorithm} \label{sec3:proposed-algorithm}
In case of a known control policy $\pi$ for a candidate DTCBF $h$ with a $\gamma \in \mathcal{K}^{\leqslant \mathrm{id}}_{\infty}$, and given that the optimization problem \eqref{eq:sec3:prop:verification-optimization-known} is generally non-convex, our algorithm is based on constructing convex underestimators of the objective function \eqref{eq:sec3:prop:verification-optimization-known-objective} and the constraint \eqref{eq:sec3:prop:verification-optimization-known-constraint} on the set $\mathbb{X}$. If its global minimum is non-negative, $(h, \gamma)$ with its friend $\pi$ is verified as a valid DTCBF. Otherwise, we partition $\mathbb{X}$ into smaller subdomains to construct tighter convex underestimators. We repeat the same procedure within each subdomain. We discard subdomains in which $(h, \gamma)$ with its friend $\pi$ is verified as a valid DTCBF or those that are entirely outside the zero-superlevel set $\mathcal{C}$ of $h$, and we further subdivide the remaining subdomains. We continue this process until either all subdomains are discarded, showing that $(h, \gamma)$ with its friend $\pi$ is a valid DTCBF, or we find a specific $x \in \mathcal{C}$ for which the DTCBF constraint \eqref{eq:sec2:DTCBF-constraint} cannot be satisfied, providing a counterexample that either $(h, \gamma)$ is not a DTCBF or $\pi$ is not a friend of $(h, \gamma)$.

If a corresponding control policy is unknown, instead of solving the min-max problem \eqref{eq:sec3:prop:verification-optimization-unknown}, we propose an alternative algorithm comprising three steps to verify whether the global optimum $\mathcal{F}^{*}$ of \eqref{eq:sec3:prop:verification-optimization-unknown} is non-negative. The essence of the min-max problem lies in finding an input $u \in \mathbb{U}$ that satisfies the DTCBF constraint \eqref{eq:sec2:DTCBF-constraint} for each \mbox{$x \in \mathcal{C}$}. In our proposed approach, we determine a suitable control input that satisfies the DTCBF constraint \eqref{eq:sec2:DTCBF-constraint} for a specific state within $\mathbb{X}$ by formulating the inner maximization in \eqref{eq:sec3:prop:verification-optimization-unknown}. If this control input satisfies the DTCBF constraint \eqref{eq:sec2:DTCBF-constraint} for all states in $\mathbb{X}$, $(h, \gamma)$ is verified as a valid DTCBF. 
Otherwise, we proceed to partition $\mathbb{X}$ into smaller subdomains and repeat the same procedure within each subdomain.
We continue this process until $(h, \gamma)$ is either verified as a valid DTCBF within each subdomain that intersects $\mathcal{C}$ and then a piecewise constant friend $\pi$ is obtained, or a counterexample is found showing that $(h, \gamma)$ is not a DTCBF. 

To simplify the introduction of the algorithm, we initially focus on the verification of a candidate DTCBF within \mbox{\textit{a subdomain}} $X^{(k)} \defeq [x^{lb, (k)}, ~ x^{ub, (k)}] \subseteq \mathbb{X}$, where $k \in \mathbb{N}_0$ is the subdomain number. 
\subsubsection{Known Control Policy -- Verification within $X^{(k)}$}
As discussed in Proposition \ref{proposition:sec3:verification-known}, for the system \eqref{eq:sec2:discrete-time-dynamical-system} with $\mathbb{U}$, $(h, \gamma)$ with a control policy $\pi:\mathcal{C} \rightarrow \mathbb{U}$ is a valid DTCBF within $X^{(k)}$, if $\mathcal{F}^{*,(k)} \geqslant 0$, where
\begin{subequations} \label{eq:sec3:verification:known-general}
    \begin{align} 
    \mathcal{F}^{*,(k)} \hspace{-0.002cm} \defeq \hspace{-0.002cm} \min_{x \in X^{(k)}} ~&~ \underbrace{h(f(x, \pi(x))) - h(x) + \gamma(h(x))}_{\mathcal{F}(x) \defeq} \label{eq:sec3:verification:known-general-objective}\\
    \mathst ~&~ \mathcal{H}(x) \defeq -h(x) \leqslant 0. \label{eq:sec3:verification:known-general-constraint}
    \end{align}
\end{subequations}
Since \eqref{eq:sec3:verification:known-general} is generally non-convex, we construct a convex underestimator $\widebreve{\mathcal{F}}^{(k)}:\mathbb{R}^n \rightarrow \mathbb{R}$ of the objective function \eqref{eq:sec3:verification:known-general-objective} and a convex underestimator $\widebreve{\mathcal{H}}^{(k)}:\mathbb{R}^n \rightarrow \mathbb{R}$ of the constraint \eqref{eq:sec3:verification:known-general-constraint} on $X^{(k)}$ as
\begin{align} 
    \widebreve{\mathcal{F}}^{(k)}(x) &\hspace{-0.05cm} \defeq \hspace{-0.05cm} \mathcal{F}(x) \hspace{-0.05cm} + \hspace{-0.05cm} \sum_{i=1}^{n} \hspace{-0.05cm} \alpha_{\mathcal{F},i}^{(k)}\bigl(x_i^{lb,(k)} \hspace{-0.1cm} - x_i\bigr)\bigl(x_i^{ub,(k)} \hspace{-0.1cm} - x_i\bigr), \label{eq:sec3:verification:known-objective-convexified} \\
    \widebreve{\mathcal{H}}^{(k)}(x) &\hspace{-0.05cm}\defeq \hspace{-0.05cm} \mathcal{H}(x) \hspace{-0.05cm} + \hspace{-0.05cm}  \sum_{i=1}^{n} \hspace{-0.05cm} \alpha_{\mathcal{H},i}^{(k)}\bigl(x_i^{lb,(k)}\hspace{-0.1cm} - x_i\bigr)\bigl(x_i^{ub,(k)} \hspace{-0.1cm} - x_i\bigr),\label{eq:sec3:verification:known-constraint-convexified}
\end{align}
where $\alpha_{\mathcal{F},i}^{(k)}, \alpha_{\mathcal{H},i}^{(k)} \hspace{-0.05cm} \in \hspace{-0.05cm} \mathbb{R}_{\geqslant 0}$, $i \hspace{-0.05cm} \in \hspace{-0.05cm} \{1, \hdots, n\}$, are sufficiently large. 
Thus, the convex optimization problem is constructed as 
\begin{subequations} \label{eq:sec4:verification:known-convex}
    \begin{align} 
    \widebreve{\mathcal{F}}^{*, (k)} \defeq \min_{x \in X^{(k)}} ~&~ \widebreve{\mathcal{F}}^{(k)}(x) \label{eq:sec3:verification:known-convex-objective}\\
    \mathst ~&~ \widebreve{\mathcal{H}}^{(k)}(x) \leqslant 0. \label{eq:sec3:verification:known-convex-constraint}
    \end{align}
\end{subequations}
Given that for all $x \in X^{(k)}$,
\begin{align}
    \widebreve{\mathcal{F}}^{(k)}(x) &\leqslant \mathcal{F}^{(k)}(x), \label{eq:verification:known:constraint-objective} \\
    \widebreve{\mathcal{H}}^{(k)}(x) &\leqslant \mathcal{H}^{(k)}(x), \label{eq:verification:known:constraint-constraint}
\end{align}
it follows that \mbox{$\widebreve{\mathcal{F}}^{*, (k)} \leqslant \mathcal{F}^{*, (k)}$}. Thus, if \eqref{eq:sec4:verification:known-convex} is feasible and $\widebreve{\mathcal{F}}^{*, (k)} \geqslant 0$, $(h,\gamma)$ with its friend $\pi$ is verified as a valid DTCBF within $X^{(k)}$ since \mbox{$\mathcal{F}^{*, (k)} \geqslant \widebreve{\mathcal{F}}^{*, (k)} \geqslant 0$}. Otherwise, we follow the steps of the algorithm discussed in ``3) \nameref{sec:sec4:algorithm-overview}'' below.
\subsubsection{Unknown Control Policy -- Verification within $X^{(k)}$} \label{sec:sec3:verification-within-domain-unknown}
For the unknown control policy case, the three-step approach to verify a candidate DTCBF $h$ with a \mbox{$\gamma \in \mathcal{K}^{\leqslant \mathrm{id}}_{\infty}$} within subdomain $X^{(k)} \defeq [x^{lb, (k)}\hspace{-0.035cm}, ~ x^{ub, (k)}] \subseteq \mathbb{X}$ for the \mbox{system \eqref{eq:sec2:discrete-time-dynamical-system}} with the control admissible set $\mathbb{U}$ is proposed as follows:
\begin{enumerate}[label=\underline{Step \Roman*}:,ref=Step \Roman*, wide=0pt] 
\item In this step, we select the state $\bar{x}^{(k)} \in \mathbb{R}^n$ as the middle of $X^{(k)}$, i.e., 
\begin{align} \label{eq:sec3:verification:unknown:step1}
    \bar{x}^{(k)}_i \defeq \frac{x^{ub, (k)}_i + x^{lb, (k)}_i}{2}, 
\end{align}
where $i\in \{1, \hdots, n\}$ and $\bar{x}^{(k)} \defeq [\bar{x}^{(k)}_1 \hspace{-0.1cm} \hdots~ \bar{x}^{(k)}_n]^{\T}$. 
\item \label{algorithm:verification:unknown:step2} For the state $\bar{x}^{(k)} \in \mathbb{R}^n$, we aim to find an admissible control input $u^{*,(k)} \in \mathbb{U}$ that satisfies the DTCBF constraint \eqref{eq:sec2:DTCBF-constraint}. To this end, we formulate the inner maximization problem in \eqref{eq:sec3:prop:verification-optimization-unknown} and solve it to global optimality using the $\alpha$BB algorithm \cite{Floudas2010a}: 
\begin{align} \label{eq:sec3:verification:unknown:step2}
    u^{*,(k)} \in {\argmax_{u \in \mathbb{U}}} \hspace{0.1cm}& \underbrace{h\bigl(f(\bar{x}^{(k)},u)\bigr) - h\bigl(\bar{x}^{(k)}\bigr) + \gamma\bigl(h(\bar{x}^{(k)})\bigr)}_{\mathcal{F}^{(k)}_u(u) \defeq}. 
\end{align}
\item \label{algorithm:verification:unknown:step3} We check whether the control input $u^{*,(k)} \in \mathbb{U}$ satisfies the DTCBF constraint \eqref{eq:sec2:DTCBF-constraint} for all $x \in X^{(k)}$. For this purpose, we check whether $\mathcal{F}^{*,(k)} \geqslant 0$, where 
\begin{subequations} \label{eq:sec3:verification:unknown:step3-general}
    \begin{align} 
        \mathcal{F}^{*,(k)}  \defeq \min_{x \in X^{(k)}}   ~& \underbrace{h\bigl(f(x,u^{*,(k)})\bigr) - h(x) +\gamma(h(x))}_{\mathcal{F}(x) \defeq} \label{eq:sec3:verification:unknown:step3-general:objective}\\
        \mathst ~& \mathcal{H}(x) \defeq -h(x) \leqslant 0.  \label{eq:sec3:verification:unknown:step3-general:constraint}
    \end{align}
\end{subequations} 
Since \eqref{eq:sec3:verification:unknown:step3-general} is generally non-convex, we construct the convex optimization problem as 
\begin{subequations}\label{eq:sec3:verification:unknown:step3-convex}
    \begin{align} 
        \widebreve{\mathcal{F}}^{*, (k)}  \defeq \min_{x \in X^{(k)}}  ~&~ \widebreve{\mathcal{F}}^{(k)}(x) \label{eq:sec3:verification:unknown:step3-convex-objective} \\
        \mathst ~&~ \widebreve{\mathcal{H}}^{(k)}(x) \leqslant 0, \label{eq:sec3:verification:unknown:step3-convex-constraint}
    \end{align}
\end{subequations}
where $\widebreve{\mathcal{F}}^{(k)}$ and $\widebreve{\mathcal{H}}^{(k)}$ are convex underestimators of the objective function \eqref{eq:sec3:verification:unknown:step3-general:objective} and the constraint \eqref{eq:sec3:verification:unknown:step3-general:constraint} on $X^{(k)}$ constructed as \eqref{eq:sec3:verification:known-objective-convexified} and \eqref{eq:sec3:verification:known-constraint-convexified}, respectively.
Similar to the known control policy case, if \eqref{eq:sec3:verification:unknown:step3-convex} is feasible and $\widebreve{\mathcal{F}}^{*, (k)} \geqslant 0$, $(h,\gamma)$ is verified as a valid DTCBF within $X^{(k)}$. Otherwise, we follow the steps of the algorithm below.
\end{enumerate}
\subsubsection{Algorithm Overview} \label{sec:sec4:algorithm-overview}
By solving the optimization problem \eqref{eq:sec4:verification:known-convex} for the known control policy case or \eqref{eq:sec3:verification:unknown:step3-convex} for the unknown control policy case on a subdomain $X^{(k)}$, we encounter three cases: 
\begin{enumerate}[label=(\Alph*), ref=(\Alph*), leftmargin=18pt]
    \item  $\widebreve{\mathcal{F}}^{*, (k)} \geqslant 0$: $(h,\gamma)$ is verified as a valid DTCBF within $X^{(k)}$ as \mbox{$\mathcal{F}^{*, (k)} \geqslant \widebreve{\mathcal{F}}^{*, (k)} \geqslant 0$}. Thus, we do not have to consider this subdomain any further. \label{case:convex-verification-known-A}
    \item The convex optimization problem \eqref{eq:sec4:verification:known-convex} or \eqref{eq:sec3:verification:unknown:step3-convex} is infeasible: Given that \eqref{eq:verification:known:constraint-constraint} holds, 
    $X^{(k)} \hspace{-0.06cm} \cap \hspace{0.05cm} \mathcal{C} = \emptyset$. Thus, we discard this subdomain. \label{case:convex-verification-known-B}
    \item \label{case:convex-verification-known-C} $\widebreve{\mathcal{F}}^{*, (k)} < 0$: For the known control policy case, consider $x^{*,(k)}$ as the global minimizer of \eqref{eq:sec4:verification:known-convex} and $\mathcal{F}$ as defined in \eqref{eq:sec3:verification:known-general}.
    For the unknown control policy case, consider $\bar{x}^{(k)}$ as obtained from \eqref{eq:sec3:verification:unknown:step1}, $\mathcal{F}^{(k)}_u$ as defined in \eqref{eq:sec3:verification:unknown:step2}, and \mbox{$u^{*,(k)}$} as a global minimizer of \eqref{eq:sec3:verification:unknown:step2}.
    \begin{enumerate} [label=(C.\arabic*), ref=(C.\arabic*), leftmargin=18pt]
        \item $x^{*,(k)} \in \mathcal{C}$ and $\mathcal{F}(x^{*,(k)}) < 0$ for the known control policy case, or $\bar{x}^{(k)} \in \mathcal{C}$, $u^{*,(k)} \in \mathbb{U}$, and $\mathcal{F}^{(k)}_u(u^{*,(k)}) < 0$ for the unknown control policy case: We terminate the algorithm and report $x^{*,(k)}$ or $\bar{x}^{(k)}$ as a counterexample. \label{case:convex-verification-known-C-1}
        \item Otherwise: $X^{(k)}$ is divided into smaller subdomains. \label{case:convex-verification-known-C-2}
    \end{enumerate}
\end{enumerate}
\vspace{-0.2cm}

Additionally, to ensure that the proposed algorithm terminates in a finite number of iterations, we impose certain stopping criteria. In case of a known control policy, we impose stopping criteria on the maximum separation between the objective function \eqref{eq:sec3:verification:known-general-objective} and its convex underestimator \eqref{eq:sec3:verification:known-convex-objective}, as well as the constraint \eqref{eq:sec3:verification:known-general-constraint} and its convex underestimator \eqref{eq:sec3:verification:known-convex-constraint}. Specifically, the algorithm terminates on subdomain $X^{(k)}$, if 
\begin{align}
    \frac{\max_i\alpha_{\mathcal{F},i}^{(k)}}{4}\sum_{i=1}^{n} \bigl(x_i^{ub,(k)} - x_i^{lb,(k)}\bigr)^2 &\leqslant \epsilon_f \label{eq:sec4:verification:stopping-criteria-1},\\
    \frac{\max_i\alpha_{\mathcal{H},i}^{(k)}}{4}\sum_{i=1}^{n}\bigl(x_i^{ub,(k)} - x_i^{lb,(k)}\bigr)^2 &\leqslant \epsilon_h \label{eq:sec4:verification:stopping-criteria-2}, 
\end{align}
where $\epsilon_f, \epsilon_h \in \mathbb{R}_{>0}$ are predefined values. For the case where a corresponding control policy is unknown, we impose the same stopping criteria, but additionally consider the size of subdomains, which determines the accuracy of the resulting control policy. In particular, the algorithm terminates on $X^{(k)}$, if \eqref{eq:sec4:verification:stopping-criteria-1}, \eqref{eq:sec4:verification:stopping-criteria-2}, and 
\begin{align}
    \delta_k^2 \defeq \sum_{i=1}^{n}\left(x_i^{ub,(k)} - x_i^{lb,(k)}\right)^2 &\leqslant \epsilon_d \label{eq:sec4:verification:stopping-criteria-3} 
\end{align}
are met, where $\delta_k \in \mathbb{R}_{>0}$ is the diagonal size of $X^{(k)}$ and $\epsilon_d \in \mathbb{R}_{>0}$ is a predefined value. \\
The proposed algorithm is detailed in Algorithm \ref{algorithm:verification}. \\
In this algorithm, if a control policy is not provided, a piecewise constant friend $\pi$ is computed. This control policy can apply a distinct control input $u \in \mathbb{U}$ for each subdomain $X^{(k)} \subseteq \mathbb{X}$ whose diagonal size $\delta_k$ satisfies $\delta_k > \sqrt{\epsilon_d}$\hspace{0.05cm}. 
\vspace{-0.1cm}
\begin{algorithm}[H]
\caption{Algorithm for Verification} \label{algorithm:verification}
\begin{algorithmic}[1] 
\State \textbf{Input:} $f$, $h$, $\gamma$, $\pi$ (if given). $\epsilon_f, \epsilon_h \in \mathbb{R}_{>0}$, and additionally $\epsilon_d \in \mathbb{R}_{>0}$ if the control policy is unknown. An $n$-rectangle set $\mathbb{X}$ such that $\mathcal{C}\subseteq \mathbb{X}$.
\State $X^{(0)} \gets \mathbb{X}$
\State $\mathcal{L} \gets \{0\}$ \Comment{$\mathcal{L}$ is the list of subdomain numbers that remain to be handled}
\State $n_{dom} \gets 0$ \Comment{Number of subdomains}
\While{$\mathcal{L} \neq \emptyset$} \label{algorithm:while}
\State $k \gets \text{any subdomain number from the list } \mathcal{L}$ 
\If{control policy is known}
    \State solve the optimization problem \eqref{eq:sec4:verification:known-convex} on $X^{(k)}$
\ElsIf{control policy is unknown}
        \State do the three steps in Section \ref{sec:sec3:verification-within-domain-unknown} for $X^{(k)}$
\EndIf
\If{Case \ref{case:convex-verification-known-A} or Case \ref{case:convex-verification-known-B}}
    \State $\mathcal{L} \gets \mathcal{L}\backslash \{k\}$ 
    \State \Goto{algorithm:while} 
\ElsIf{Case \ref{case:convex-verification-known-C-1}}
    \State \textbf{print} ``$(h,\gamma)$ is not a valid DTCBF {\small$[$or $\pi$, if given, is not a friend$]$}, and the counterexample is:" $x^{*,(k)}$ {\small$[$or $\bar{x}^{(k)}]$}
    \State terminate the algorithm
\Else \Comment{Case \ref{case:convex-verification-known-C-2}}
\If{stopping criteria are met on $X^{(k)}$}
    \State \textbf{print} ``$(h,\gamma)$ is not valid {\small$[$or $\pi$, if given, is not a friend$]$} for the selected conservatism $\epsilon_f, \epsilon_h$ {\small$[$and $\epsilon_d]$}"
    \State terminate the algorithm
\EndIf
\State divide $X^{(k)}$ into $X^{(n_{dom}+1)}$ and $X^{(n_{dom}+2)}$
\State $\mathcal{L} \gets (\mathcal{L}\backslash\{k\}) \cup \{n_{dom}+1,~ n_{dom}+2\}$
\State $n_{dom} \gets n_{dom}+2$
\EndIf
\EndWhile
\State \textbf{print} ``$(h,\gamma)$ is a DTCBF"
\end{algorithmic}
\end{algorithm}
\vspace{-0.1cm}
\begin{remark} \label{remark:differences-ourmethod-abb}
The main differences between our method and the $\alpha$BB algorithm \cite{Floudas2010a} are as follows:
    \begin{itemize}[leftmargin=6.4pt]
        \item For the known control policy case, we only compute lower bounds on the global minimum of \eqref{eq:sec3:prop:verification-optimization-known} (through optimization problem \eqref{eq:sec4:verification:known-convex}), not upper bounds. We also discard subdomain $X^{(k)}$, where $\widebreve{\mathcal{F}}^{*, (k)} \geqslant 0$. Thus, our method can be significantly faster compared to the $\alpha$BB algorithm for the verification problem. Using our method, we can falsify candidate DTCBFs by providing counterexamples without any conservatism, which is not possible using the $\alpha$BB algorithm due to the concept of $\epsilon$-feasibility (see \mbox{Section \ref{sec:sec5}}). 
        \item For the unknown control policy case, the $\alpha$BB algorithm is incapable of solving \eqref{eq:sec3:prop:verification-optimization-unknown} and verifying or falsifying a candidate DTCBF.
    \end{itemize}
\end{remark}
\begin{remark} 
The computation time of Algorithm \ref{algorithm:verification} increases exponentially with each increment in the dimension of the system \eqref{eq:sec2:discrete-time-dynamical-system}. Thus, it may not be suitable for high-dimensional systems.
\end{remark}

\section{Numerical Case Study} \label{sec:sec5}
\begin{figure*}[ht!] 
\begin{subfigure}{.33\textwidth}
\centering
\includegraphics[width=1\textwidth]{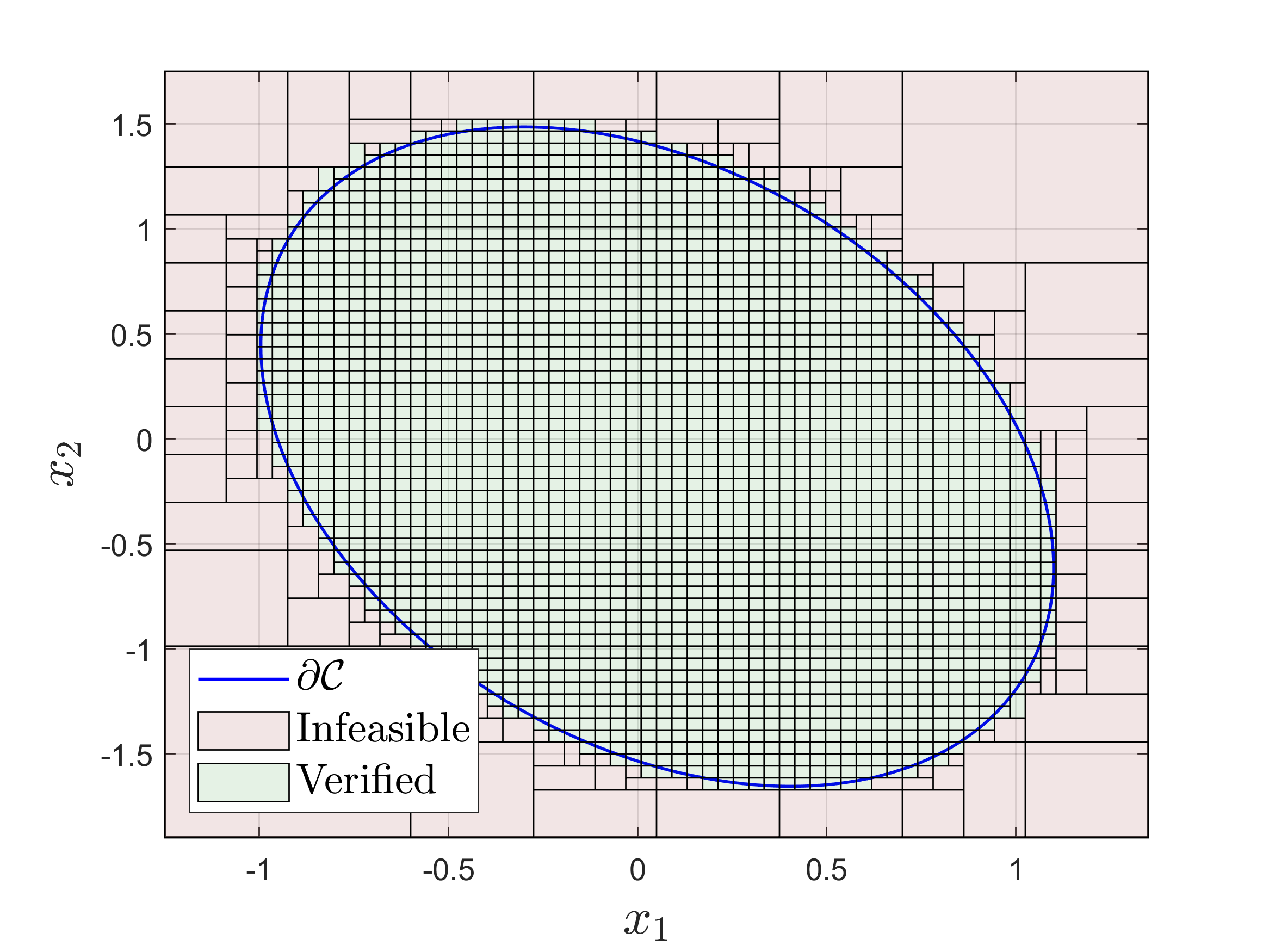} \caption{} \label{fig:example-DTCBF:1a}
\end{subfigure}
\begin{subfigure}{.33\textwidth}
\centering
\includegraphics[width=1\textwidth]{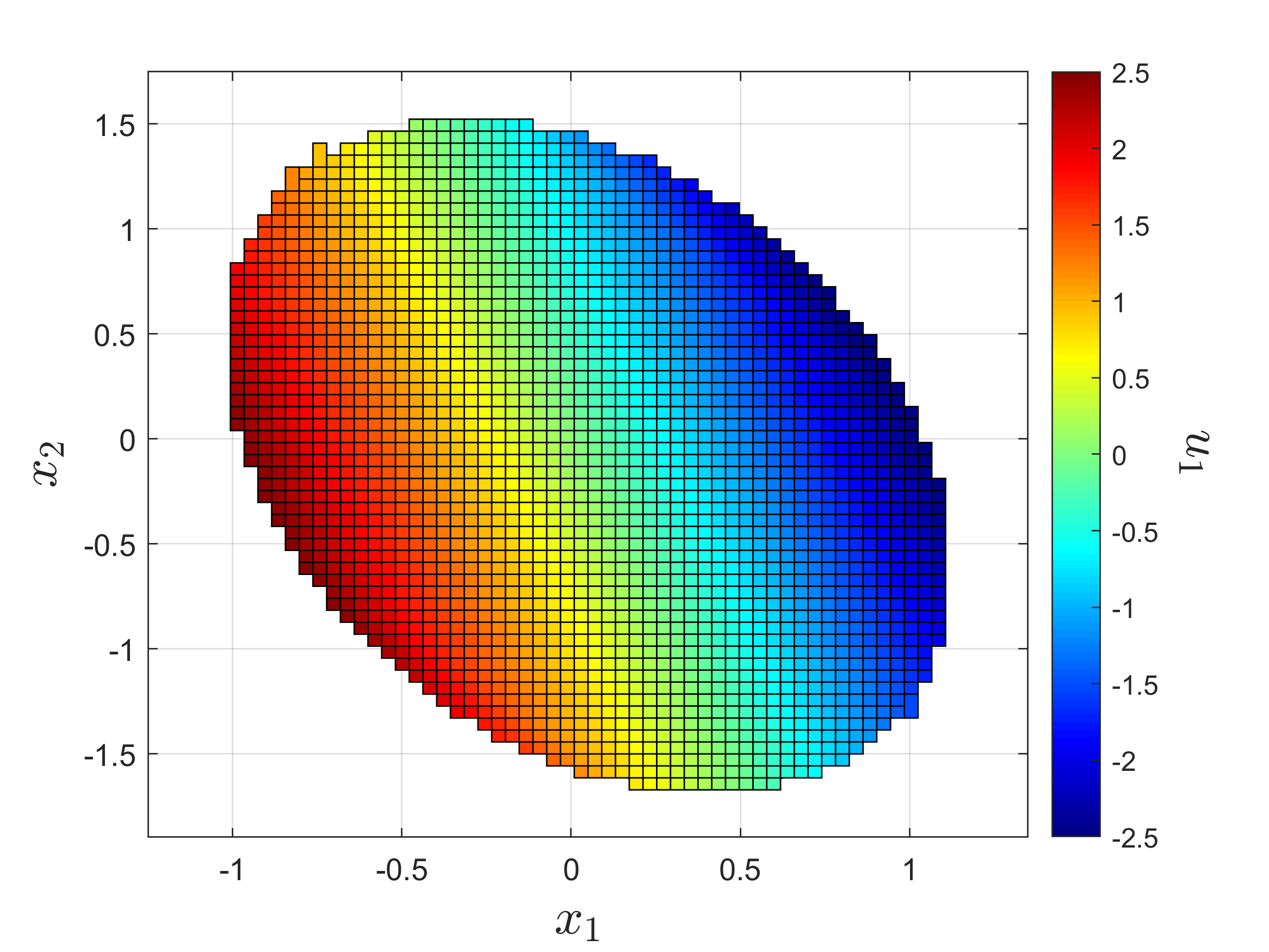} \caption{} \label{fig:example-DTCBF:1b}
\end{subfigure}
\begin{subfigure}{.33\textwidth}
\centering
\includegraphics[width=1\textwidth]{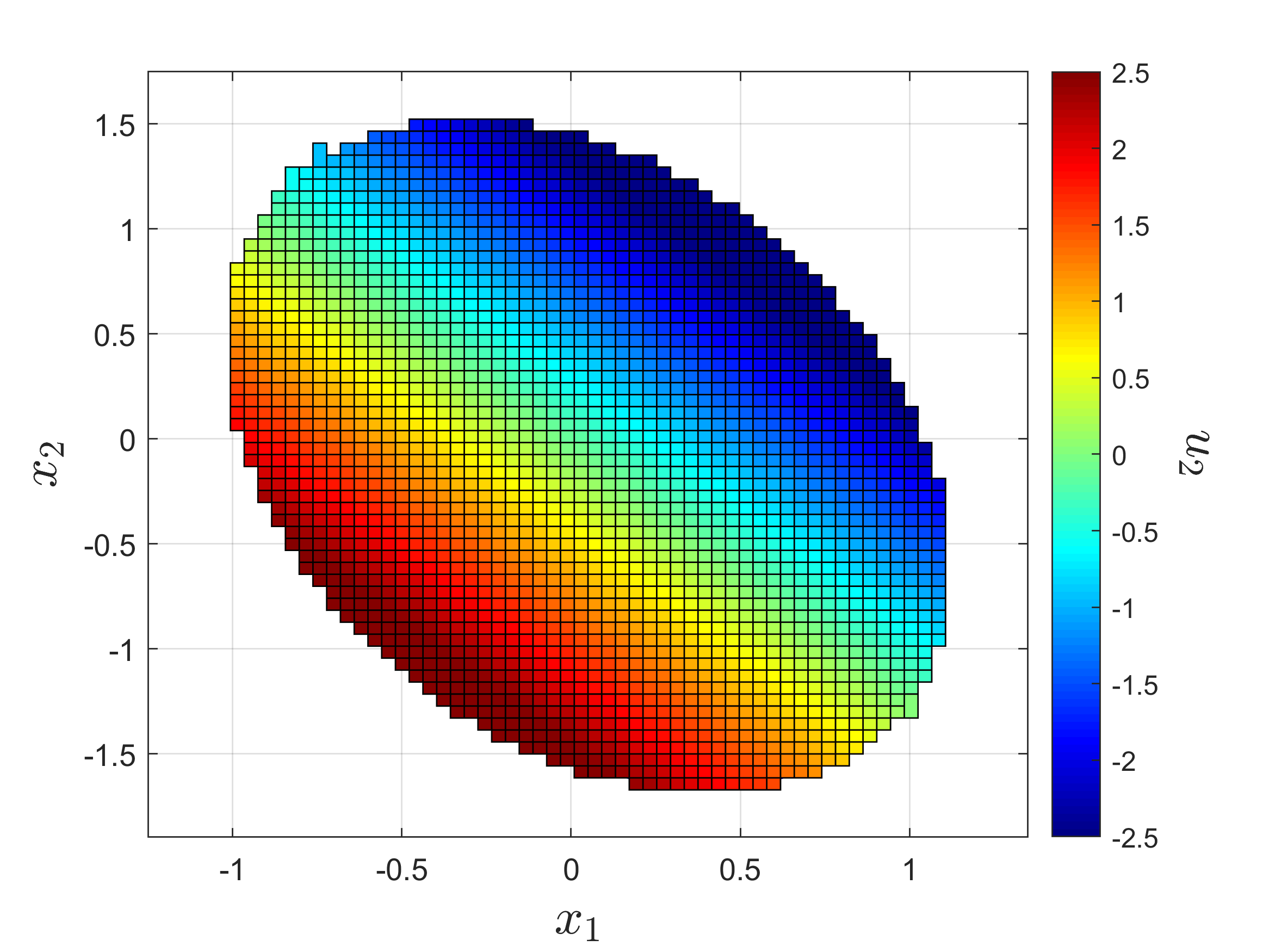} \caption{} \label{fig:example-DTCBF:1c}
\end{subfigure}
\caption{Algorithm \ref{algorithm:verification} is applied to the candidate DTCBF \eqref{eq:sec5:example:DTCBF} with $\gamma$, defined as $\gamma(r) \defeq 0.8r$, $r \in \mathbb{R}_{\geqslant 0}$, for the discrete-time system \eqref{eq:sec5:example:discrete-dynamical}. \mbox{(a) The green rectangles} signify that $(h,\gamma)$ is verified within the corresponding subdomain, Case \ref{case:convex-verification-known-A}, while the red rectangles indicate that the optimization problem \eqref{eq:sec3:verification:unknown:step3-convex} is infeasible, Case \ref{case:convex-verification-known-B}, in the sense that the corresponding subdomain is entirely outside $\mathcal{C}$. The blue curve represents the boundary of $\mathcal{C}$. (b, c) The colored subdomains represent the computed values of control inputs on each subdomain, that is, $u_1$ in (b) and $u_2$ in (c).} \vspace{-0.4cm}
\label{fig:sec5:verification-example}
\end{figure*}
In this section, we aim to compare the standard $\alpha$BB algorithm \cite{Floudas2010a} and our proposed methods for the verification problem. To this end, we revisit the continuous-time system discussed in \cite{WANG2023}:
\begin{align} \label{eq:sec5:example:continuous-dynamical}
    \begin{bmatrix}
        \dot{x}_1 \\ \dot{x}_2
    \end{bmatrix} = \begin{bmatrix}
        2 & 1 \\
        3 & 1
    \end{bmatrix}\begin{bmatrix}
        x_1 \\ x_2
    \end{bmatrix} + \begin{bmatrix}
        u_1 \\ u_2
    \end{bmatrix},
\end{align}
where $u_1 \in [-2.5,~ 2.5]$ and $u_2 \in [-2.5,~ 2.5]$. 
In \cite{WANG2023}, the continuous-time CBF is synthesized as 
\begin{align} \label{eq:sec5:example:DTCBF}
    h(x) &= −7.635x_1^2 - 3.439x_1x_2 − 3.4024x_2^2 +  \nonumber \\
    &\qquad \qquad 0.5x_1 − 0.4x_2 + 7.402,
\end{align}
and its corresponding control policy as 
\begin{subequations} \label{eq:sec5:example:Control-Policy}
\begin{align}
    \pi_1(x) &= −2.32x_1 − 1.11x_2 + 0.022,\\
    \pi_2(x) & = −2.12x_1 − 1.27x_2 − 0.046,
\end{align}
\end{subequations}
with $\pi(x) \defeq [\pi_1(x)~ \pi_2(x)]^{\T}$.
Consider the discrete-time representation of \eqref{eq:sec5:example:continuous-dynamical} using exact zero-order hold (ZOH) discretization with sample time $T_s \defeq \unit[1]{s}$, leading to 
\begin{align} \label{eq:sec5:example:discrete-dynamical}
    \begin{bmatrix}
        x^+_{1} \\ x^+_{2}
    \end{bmatrix} = \begin{bmatrix}
        17.6 & 7.3 \\
        22.0 & 10.3
    \end{bmatrix}\begin{bmatrix}
        x_{1} \\ x_{2}
    \end{bmatrix} + \begin{bmatrix}
        5.4 & 2.0 \\
        5.9 & 3.4
    \end{bmatrix} \begin{bmatrix}
        u_{1} \\ u_{2}
    \end{bmatrix}.
\end{align}
We aim to verify whether $(h, \gamma)$, where $\gamma(r) \defeq 0.8r$, \mbox{$r \in \mathbb{R}_{\geqslant 0}$}, with the control policy $\pi$ as in \eqref{eq:sec5:example:Control-Policy}, is a valid DTCBF for the discrete-time system \eqref{eq:sec5:example:discrete-dynamical}.

Using the $\alpha$BB algorithm for the convergence tolerance \mbox{$\epsilon_c \defeq 10^{-6}$} and the feasibility tolerance $\epsilon_f \defeq 10^{-12}$ (both related to the stopping criteria of the $\alpha$BB algorithm as detailed in \cite{Floudas2010a}), a global minimizer $x^*$ of the optimization problem \eqref{eq:sec3:prop:verification-optimization-known} is obtained after 55 iterations, with each iteration taking an average time of $\unit[0.4]{s}$, as 
\begin{align}
    x^*_1 = 0.841, \quad x^*_2 = -1.457. \nonumber
\end{align} 
It is evident that while the DTCBF constraint \eqref{eq:sec2:DTCBF-constraint} is not satisfied at $x^*$, where $x^* \defeq [x^*_1 ~ x^*_2]^{\T}$, the state $x^*$ lies slightly outside the zero-superlevel set $\mathcal{C}$ of $h$, in the sense that $h(x^*) < 0$. Thus, drawing a conclusion on the validity of $(h,\gamma)$ with the control policy $\pi$ remains inconclusive using the $\alpha$BB algorithm (see Remark \ref{remark:differences-ourmethod-abb}).

Then, we apply Algorithm \ref{algorithm:verification} for the known control policy case with \mbox{$\epsilon_f \defeq \epsilon_h \defeq 10^{-6}$}, using the \textit{scaled longest side} \cite{Maranas1994b} as a branching strategy. After 2 iterations, with each iteration taking an average of $\unit[0.2]{s}$, 
the output of our method presents a counterexample given by
\begin{align}
\bar{x}_1 = 1.030, \quad \bar{x}_2 = -1.110. \nonumber
\end{align}
It can be observed that the DTCBF constraint \eqref{eq:sec2:DTCBF-constraint} is not satisfied at $\bar{x}$, where $\bar{x} \defeq [\bar{x}_1 ~ \bar{x}_2]^{\T}$, and $h(\bar{x}) \geqslant 0$. Thus, either $(h,\gamma)$ is not a valid DTCBF for the system \eqref{eq:sec5:example:discrete-dynamical}, or $\pi$ is not a friend of $(h,\gamma)$.

To determine whether $(h,\gamma)$ is a valid DTCBF with a different control policy, we apply Algorithm \ref{algorithm:verification} for the unknown control policy case with \mbox{$\epsilon_f \defeq \epsilon_h \defeq \epsilon_d \defeq 10^{-6}$}, using the scaled longest side as a branching strategy. 
Application of the proposed algorithm results in Figure \ref{fig:sec5:verification-example}. As observed in Figure \ref{fig:example-DTCBF:1a}, either the subdomains are outside of $\mathcal{C}$, or $(h,\gamma)$ is verified within the subdomains. As a result, $(h,\gamma)$ with the computed piecewise constant friend, as depicted in Figures \ref{fig:example-DTCBF:1b} and \ref{fig:example-DTCBF:1c}, is verified as a valid DTCBF after approximately 2500 iterations, each taking an average of $\unit[0.4]{s}$. Simulations are run on an Intel i7-7700HQ machine with Matlab 2022a.


\section{Conclusions and Future Work} \label{sec:sec6}
In this paper, we have proposed a novel branch-and-bound method inspired by the 
$\alpha$BB algorithm for the verification of candidate DTCBFs. This method either verifies a candidate function as a valid DTCBF for a discrete-time system with input constraints or falsifies it by providing a counterexample. The method is applicable in both cases, whether a corresponding control policy is known or unknown. We have applied our method to a numerical case study. 

In our future work, we aim to extend this method to verify candidate \textit{robust} DTCBFs, thereby guaranteeing robust safety for dynamical systems under disturbances. Additionally, we aim to propose a method to \textit{synthesize} DTCBFs, where we envision that the current verification algorithm and the obtained insights are instrumental.
\vspace{-0.03cm}

\bibliographystyle{IEEEtran}
\bibliography{IEEEabrv}

\begin{thebibliography}{10}
\providecommand{\url}[1]{#1}
\csname url@rmstyle\endcsname
\providecommand{\newblock}{\relax}
\providecommand{\bibinfo}[2]{#2}
\providecommand\BIBentrySTDinterwordspacing{\spaceskip=0pt\relax}
\providecommand\BIBentryALTinterwordstretchfactor{4}
\providecommand\BIBentryALTinterwordspacing{\spaceskip=\fontdimen2\font plus
\BIBentryALTinterwordstretchfactor\fontdimen3\font minus \fontdimen4\font\relax}
\providecommand\BIBforeignlanguage[2]{{%
\expandafter\ifx\csname l@#1\endcsname\relax
\typeout{** WARNING: IEEEtran.bst: No hyphenation pattern has been}%
\typeout{** loaded for the language `#1'. Using the pattern for}%
\typeout{** the default language instead.}%
\else
\language=\csname l@#1\endcsname
\fi
#2}}

\bibitem{Ames2014a}
A.~D. Ames, J.~W. Grizzle, and P.~Tabuada, ``Control barrier function based quadratic programs with application to adaptive cruise control,'' in \emph{IEEE Conference on Decision and Control}, 2014, pp. 6271--6278.

\bibitem{Agrawal2017a}
A.~Agrawal and K.~Sreenath, ``{Discrete Control Barrier Functions for Safety-Critical Control of Discrete Systems with Application to Bipedal Robot Navigation},'' in \emph{Robotics: Science and Systems}, 2017.

\bibitem{Zeng2021a}
J.~Zeng, B.~Zhang, and K.~Sreenath, ``{Safety-Critical Model Predictive Control with Discrete-Time Control Barrier Function},'' in \emph{IEEE American Control Conference}, 2021, pp. 3882--3889.

\bibitem{Zeng2021b}
J.~Zeng, Z.~Li, and K.~Sreenath, ``{Enhancing Feasibility and Safety of Nonlinear Model Predictive Control with Discrete-Time Control Barrier Functions},'' in \emph{IEEE Conference on Decision and Control}, 2021, pp. 6137--6144.

\bibitem{katriniok2023}
A.~Katriniok, E.~Shakhesi, and W.~Heemels, ``Discrete-time control barrier functions for guaranteed recursive feasibility in nonlinear mpc: An application to lane merging,'' in \emph{2023 62nd IEEE Conference on Decision and Control (CDC)}, 2023, pp. 3776--3783.

\bibitem{Didier2023}
A.~Didier, R.~C. Jacobs, J.~Sieber, K.~P. Wabersich, and M.~N. Zeilinger, ``Approximate predictive control barrier functions using neural networks: A computationally cheap and permissive safety filter,'' in \emph{2023 European Control Conference (ECC)}, 2023, pp. 1--7.

\bibitem{Zhang2023}
H.~Zhang, Z.~Li, H.~Dai, and A.~Clark, ``Efficient sum of squares-based verification and construction of control barrier functions by sampling on algebraic varieties,'' in \emph{2023 62nd IEEE Conference on Decision and Control (CDC)}, 2023, pp. 5384--5391.

\bibitem{Clark2021}
A.~Clark, ``Verification and synthesis of control barrier functions,'' in \emph{2021 60th IEEE Conference on Decision and Control (CDC)}, 2021, pp. 6105--6112.

\bibitem{Clark2022}
------, ``A semi-algebraic framework for verification and synthesis of control barrier functions,'' 2022, \textit{arXiv preprint arXiv:2209.00081}.

\bibitem{Isaly2022}
A.~Isaly, M.~Ghanbarpour, R.~G. Sanfelice, and W.~E. Dixon, ``On the feasibility and continuity of feedback controllers defined by multiple control barrier functions for constrained differential inclusions,'' in \emph{2022 American Control Conference (ACC)}, 2022, pp. 5160--5165.

\bibitem{Tan2022}
X.~Tan and D.~V. Dimarogonas, ``Compatibility checking of multiple control barrier functions for input constrained systems,'' in \emph{2022 IEEE 61st Conference on Decision and Control (CDC)}, 2022, pp. 939--944.

\bibitem{Floudas2010a}
C.~A. Floudas, \emph{The $\alpha$BB Approach for General Constrained Twice-Differentiable NLPs : Theory}.\hskip 1em plus 0.5em minus 0.4em\relax Boston, MA: Springer US, 2000, pp. 333--375.

\bibitem{dReal}
S.~Gao, S.~Kong, and E.~M. Clarke, ``dreal: An smt solver for nonlinear theories over the reals,'' in \emph{Automated Deduction -- CADE-24}, M.~P. Bonacina, Ed.\hskip 1em plus 0.5em minus 0.4em\relax Berlin, Heidelberg: Springer Berlin Heidelberg, 2013, pp. 208--214.

\bibitem{wu2023neural}
J.~Wu, A.~Clark, Y.~Kantaros, and Y.~Vorobeychik, ``Neural lyapunov control for discrete-time systems,'' in \emph{Advances in Neural Information Processing Systems}, vol.~36.\hskip 1em plus 0.5em minus 0.4em\relax Curran Associates, Inc., 2023, pp. 2939--2955.

\bibitem{Ames2019a}
A.~D. Ames, S.~Coogan, M.~Egerstedt, G.~Notomista, K.~Sreenath, and P.~Tabuada, ``{Control Barrier Functions: Theory and Applications},'' in \emph{European Control Conference}, 2019, pp. 3420--3431.

\bibitem{chen2021}
Y.~Chen, M.~Jankovic, M.~Santillo, and A.~D. Ames, ``Backup control barrier functions: Formulation and comparative study,'' in \emph{2021 60th IEEE Conference on Decision and Control (CDC)}, 2021, pp. 6835--6841.

\bibitem{tonkens2022refining}
S.~Tonkens and S.~Herbert, ``Refining control barrier functions through hamilton-jacobi reachability,'' in \emph{2022 IEEE/RSJ International Conference on Intelligent Robots and Systems (IROS)}, 2022, pp. 13\,355--13\,362.

\bibitem{ADJIMAN1998a}
C.~Adjiman, S.~Dallwig, C.~Floudas, and A.~Neumaier, ``A global optimization method, $\alpha$bb, for general twice-differentiable constrained nlps — i. theoretical advances,'' \emph{Computers \& Chemical Engineering}, vol.~22, no.~9, pp. 1137--1158, 1998.

\bibitem{Maranas1994b}
C.~D. Maranas and C.~A. Floudas, ``Global minimum potential energy conformations of small molecules,'' \emph{Journal of Global Optimization}, vol.~4, pp. 135--170, 1994.

\bibitem{WANG2023}
H.~Wang, K.~Margellos, and A.~Papachristodoulou, ``Safety verification and controller synthesis for systems with input constraints,'' \emph{IFAC-PapersOnLine}, vol.~56, no.~2, pp. 1698--1703, 2023, 22nd IFAC World Congress.

\end{thebibliography}

\end{document}